\RequirePackage{fix-cm}
\documentclass[smallextended]{svjour3}       % onecolumn (second format)
\smartqed  % flush right qed marks, e.g. at end of proof
\usepackage{amsmath,amssymb,amsfonts}
\usepackage[english]{babel}
\usepackage{graphicx}
\usepackage{mathtools}
\usepackage{ctable}
\usepackage{latexsym}
\usepackage{color}
\usepackage{amsfonts}
\usepackage{graphicx}
\begin{document}

\title{On cumulative entropies in terms of moments of order statistics%\thanks{Grants or other notes
%about the article that should go on the front page should be
%placed here. General acknowledgments should be placed at the end of the article.}
}
%\subtitle{Do you have a subtitle?\\ If so, write it here}

%\titlerunning{Short form of title}        % if too long for running head

\author{Narayanaswamy Balakrishnan         \and
        Francesco Buono \and Maria Longobardi %etc.
}

%\authorrunning{Short form of author list} % if too long for running head

\institute{N. Balakrishnan, \at
             McMaster University, Hamilton, Ontario, Canada \\
              \email{bala@mcmaster.ca}           %  \\
%             \emph{Present address:} of F. Author  %  if needed
           \and
           F. Buono \at
             Università degli Studi di Napoli Federico II, Naples, Italy  \\
     \email{francesco.buono3@unina.it}
           \and
          M. Longobardi \at
             Università degli Studi di Napoli Federico II, Naples, Italy  \\
     \email{maria.longobardi@unina.it}
}

\date{Received: date / Accepted: date}
% The correct dates will be entered by the editor

\maketitle

\begin{abstract}
In this paper relations among some kinds of cumulative entropies and moments of order statistics are presented. By using some characterizations and the symmetry of a non negative and absolutely continuous random variable $X$, lower and upper bounds for entropies are obtained and examples are given.  
% \PACS{PACS code1 \and PACS code2 \and more}

\end{abstract}

\section{Introduction}

In reliability theory, to describe and study the information related to a non-negative absolutely continuous random variable $X$ we use the Shannon entropy, or differential entropy, of $X$, defined by (Shannon, 1948)
\begin{equation}
\nonumber
H(X)=-\mathbb E[\log(X)]=-\int_0^{+\infty}f(x)\log f(x)\mathrm dx,
\end{equation}
where $\log$ is the natural logarithm and $f$ is the probability density function (pdf) of $X$. In the following we use $F$ and $\overline F$ to indicate the cumulative distribution function (cdf) and the survival function (sf) of $X$, respectively.

In the literature, there are several different versions of entropy, each one suitable for a specific situation. Rao et al., 2004, introduced the Cumulative Residual Entropy (CRE) of $X$ as
\begin{eqnarray}
\nonumber
\mathcal{E}(X)&=&-\int_0^{+\infty}\overline F(x)\log \overline F(x)\mathrm dx.
\end{eqnarray}

Di Crescenzo and Longobardi, 2009, introduced the Cumulative Entropy (CE) of $X$ as
\begin{eqnarray}
\label{CE}
\mathcal{CE}(X)&=&-\int_0^{+\infty}F(x)\log F(x)\mathrm dx .
\end{eqnarray}
This information measure is suitable to measure information when uncertainty is related
to the past, a dual concept of the cumulative residual entropy which relates to uncertainty on the future lifetime of a system.

Mirali et al., 2016, introduced the Weighted Cumulative Residual Entropy (WCRE) of $X$ as
\begin{eqnarray}
\nonumber
\mathcal{E}^w(X)&=&-\int_0^{+\infty}x(1-F(x))\log(1-F(x))\mathrm dx .
\end{eqnarray}

Mirali and Baratpour, 2017, introduced the Weighted Cumulative Entropy (WCE) of $X$ as
\begin{eqnarray}
\nonumber
\mathcal{CE}^w(X)&=&-\int_0^{+\infty}xF(x)\log F(x)\mathrm dx.
\end{eqnarray}

Recently, various authors have discussed different versions of entropy and their applications (see, for instance, \cite{cali}, \cite{CaLoNa}, \cite{CaLoPsa}, \cite{Lo2014}).

The paper is organized as follows. In Section 2, we study relations among some kinds of  entropies and moments of order statistics and present various examples. In Section 3, bounds are given by using also some characterizations and properties (as the symmetry) of the random variable $X$, some examples and bounds for known distributions are given.

\section{A relation among entropies and order statistics}

We recall that, if we have $n$ i.i.d. random variables $X_1,X_2,\dots,X_n$, we can introduce the order statistics $X_{k:n}$, $k=1,\dots,n$. The $k$-th order statistic is equal to the $k$-th smallest value from the sample. We know that the cdf of $X_{k:n}$ can be given in terms of the cdf of the parent distribution; in fact
\begin{equation}
\nonumber
F_{k:n}(x)=\sum_{j=k}^{n} \binom nj [ F(x) ]^{j} [ 1 - F(x) ]^{n-j},
\end{equation}
whereas the pdf of $X_{k:n}$ is expressed as
\begin{equation}
\nonumber
f_{k:n}(x)= \binom nk k[ F(x) ]^{k-1} [ 1 - F(x) ]^{n-k} f(x).
\end{equation}
Choosing $k=1$ and $k=n$ we get the smallest and the largest order statistic, respectively. Their cdf and pdf are given by
\begin{align*}
&F_{1:n}(x)=1-[1-F(x)]^n \ \ \ \ \ f_{1:n}(x)=n[1-F(x)]^{n-1}f(x) \\
&F_{n:n}(x)=[F(x)]^n \ \ \ \ \ \ \ \ \ \ \ \ \ \ \  f_{n:n}(x)=n[F(x)]^{n-1}f(x).
\end{align*}

\subsection{Cumulative residual entropy}

The Cumulative Residual Entropy (CRE) of $X$ can be written also in terms of order statistics, that is
\begin{eqnarray}
\nonumber
\mathcal{E}(X)&=&-\int_0^{+\infty}(1-F(x))\log(1-F(x))\mathrm dx \\
\nonumber
&=&-x(1-F(x))\log(1-F(x))\big|_0^{+\infty}-\int_0^{+\infty}x\log(1-F(x))f(x)\mathrm dx- \\
\nonumber
& &   \int_0^{+\infty}xf(x)\mathrm dx \\
\nonumber
&=&\int_0^{+\infty}x[-\log(1-F(x))]f(x)\mathrm dx - \mathbb E(X) \\
\nonumber
&=&\int_0^{+\infty}x\left[\sum_{n=1}^{+\infty}\frac{F(x)^n}{n}\right]f(x)\mathrm dx-\mathbb E(X) \\
\label{eq1}
&=&\sum_{n=1}^{+\infty}\frac{1}{n(n+1)}\mu_{n+1:n+1}-\mathbb E(X),
\end{eqnarray}
where $\mathbb E(X)$ is the expectation or mean of $X$, and $\mu_{n+1:n+1}$ the mean of the largest order statistic in a sample of size $n+1$ from $F$, provided that $\lim_{x\to+\infty}-x(1-F(x))\log(1-F(x))=0$.

We note that \eqref{eq1} can be rewritten as
\begin{equation}
\label{eq2}
\mathcal{E}(X)=\sum_{n=1}^{+\infty}\left(\frac{1}{n}-\frac{1}{n+1}\right)\mu_{n+1:n+1}-\mathbb E(X).
\end{equation}

\begin{example}
Consider the standard exponential distribution with pdf $f(x)=\mathrm e^{-x}$, $x>0$. Then, it is known that
\begin{equation}
\nonumber
\mathbb E(X)=1 \ \mbox{ and }\  \mathbb E(X_{n:n})=1+\frac{1}{2}+\dots+\frac{1}{n}.
\end{equation}
Then, from \eqref{eq2}, we readily have
\begin{eqnarray}
\nonumber
\mathcal E(X)&=&\left(1-\frac{1}{2}\right)\mu_{2:2}+\left(\frac{1}{2}-\frac{1}{3}\right)\mu_{3:3}+\left(\frac{1}{3}-\frac{1}{4}\right)\mu_{4:4}+\dots -\mathbb E(X) \\
\nonumber
&=& \left(\mu_{2:2}-\mathbb E(X)\right)+\frac{1}{2} \left(\mu_{3:3}-\mu_{2:2}\right)+\frac{1}{3} \left(\mu_{4:4}-\mu_{3:3}\right)+\dots \\
\nonumber
&=& \frac{1}{2}+ \frac{1}{2\cdot 3}+\frac{1}{3\cdot 4}+\dots \\
\nonumber
&=& \sum_{n=1}^{+\infty}\frac{1}{n(n+1)}=1.
\end{eqnarray}
\end{example}

\begin{example}
Consider the standard uniform distribution with pdf $f(x)=1$, $0<x<1$. Then, it is known that 
\begin{equation}
\nonumber
\mathbb E(X)=\frac{1}{2} \ \mbox{ and }\  \mathbb E(X_{n:n})=\frac{n}{n+1}.
\end{equation}
So, from \eqref{eq1}, we readily find
\begin{eqnarray}
\nonumber
\mathcal E(X)&=& \sum_{n=1}^{+\infty} \frac{1}{n(n+1)}\ \frac{n+1}{n+2}-\frac{1}{2} \\
\nonumber
&=& \frac{1}{2}\sum_{n=1}^{+\infty}\left(\frac{1}{n}-\frac{1}{n+2}\right)-\frac{1}{2} \\
\nonumber
&=& \frac{1}{2}\left(1+\frac{1}{2}\right)-\frac{1}{2}=\frac{1}{4}.
\end{eqnarray}
\end{example}

\subsection{Cumulative entropy}

The Cumulative Entropy (CE) of $X$ can be rewritten in terms of the mean of the minimum order statistic; integrating by parts \eqref{CE}
\begin{eqnarray}
\nonumber
\mathcal{CE}(X)&=&-xF(x)\log F(x)\big|_0^{+\infty}+\int_0^{+\infty}x\log F(x)f(x)\mathrm dx+\int_0^{+\infty}xf(x)\mathrm dx \\
\nonumber
&=&\int_0^{+\infty}x\log[1-(1-F(x))]f(x)\mathrm dx + \mathbb E(X) \\
\nonumber
&=&-\int_0^{+\infty}x\sum_{n=1}^{+\infty}\frac{(1-F(x))^n}{n}f(x)\mathrm dx+\mathbb E(X) \\
\label{eq3}
&=&-\sum_{n=1}^{+\infty}\frac{1}{n(n+1)}\mu_{1:n+1}+\mathbb E(X),
\end{eqnarray}
where $\mu_{1:n+1}$ is the mean of the smallest order statistic from a sample of size $n+1$ from $F$, provided that $\lim_{x\to+\infty}-xF(x)\log F(x)=0$.

We note that \eqref{eq3} can be rewritten as
\begin{equation}
\label{eq4}
\mathcal{CE}(X)=-\sum_{n=1}^{+\infty}\left(\frac{1}{n}-\frac{1}{n+1}\right)\mu_{1:n+1}+\mathbb E(X).
\end{equation}

\begin{example}
For the standard exponential distribution, it is known that
\begin{equation}
\nonumber
\mathbb \mu_{1:n}=\frac{1}{n},
\end{equation}
and so from \eqref{eq4}, we readily have
\begin{eqnarray}
\nonumber
\mathcal{CE}(X)&=&-\sum_{n=1}^{+\infty}\left(\frac{1}{n}-\frac{1}{n+1}\right)\frac{1}{n+1}+1 \\
\nonumber
&=& -\sum_{n=1}^{+\infty}\frac{1}{n(n+1)}+\sum_{n=1}^{+\infty}\frac{1}{(n+1)^2}+1 \\
\nonumber
&=& \frac{\pi^2}{6}-1,
\end{eqnarray}
by the use of Euler's identity.
\end{example}

\begin{example}
For the standard uniform distribution, using the fact that
\begin{equation}
\nonumber
\mathbb \mu_{1:n}=\frac{1}{n+1},
\end{equation}
we obtain from \eqref{eq3} that
\begin{eqnarray}
\nonumber
\mathcal{CE}(X)&=&-\sum_{n=1}^{+\infty}\frac{1}{n(n+1)(n+2)}+\frac{1}{2} \\
\nonumber
&=& \frac{1}{2}-\frac{1}{2}\sum_{n=1}^{+\infty}\frac{1}{n}+\sum_{n=1}^{+\infty}\frac{1}{n+1}-\frac{1}{2}\sum_{n=1}^{+\infty}\frac{1}{n+2} \\
\nonumber
&=& \frac{1}{4}
\end{eqnarray}
by the use of Euler's identity.
\end{example}

\begin{remark}
If the random variable $X$ has finite mean $\mu$ and is symmetrically distributed about $\mu$, then we know
\begin{equation}
\nonumber
\mu_{n:n}-\mu=\mu-\mu_{1:n},
\end{equation}
and so the symmetry property of $\mathcal{CE}$ readily follows.
\end{remark}

\subsection{Weighted cumulative entropies}

In the same way the Weighted Cumulative Residual Entropy (WCRE) of $X$ can be expressed as
\begin{eqnarray}
\nonumber
\mathcal{E}^w(X)&=&-\frac{x^2}{2}(1-F(x))\log(1-F(x))\big|_0^{+\infty}-\frac{1}{2}\int_0^{+\infty}x^2\log(1-F(x))f(x)\mathrm dx \\
\nonumber
&&-\frac{1}{2}\int_0^{+\infty}x^2f(x)\mathrm dx \\
\nonumber
&=&\frac{1}{2}\int_0^{+\infty}x^2[-\log(1-F(x))]f(x)\mathrm dx - \frac{1}{2}\mathbb E(X^2) \\
\nonumber
&=&\frac{1}{2}\int_0^{+\infty}x^2\left[\sum_{n=1}^{+\infty}\frac{F(x)^n}{n}\right]f(x)\mathrm dx - \frac{1}{2}\mathbb E(X^2) \\
\label{eq5}
&=&\frac{1}{2}\sum_{n=1}^{+\infty}\frac{1}{n(n+1)}\mu^{(2)}_{n+1:n+1} - \frac{1}{2}\mathbb E(X^2),
\end{eqnarray}
where $\mu^{(2)}_{n+1:n+1}$ is the second moment of the largest order statistic in a sample of size $n+1$, provided that $\lim_{x\to+\infty}-\frac{x^2}{2}(1-F(x))\log(1-F(x))=0$.

\begin{example}
For the standard uniform distribution, using the fact that 
\begin{equation}
\nonumber
\mu^{(2)}_{n+1:n+1}=\frac{n+1}{n+3}\  \mbox{ and } \ \mathbb E(X^2)=\frac{1}{3},
\end{equation}
we obtain from \eqref{eq5} that
\begin{eqnarray}
\nonumber
\mathcal E^w(X)&=& \frac{1}{2}\sum_{n=1}^{+\infty} \frac{1}{n(n+1)}\ \frac{n+1}{n+3}-\frac{1}{6} \\
\nonumber
&=& \frac{1}{6}\sum_{n=1}^{+\infty}\left(\frac{1}{n}-\frac{1}{n+3}\right)-\frac{1}{6} \\
\nonumber
&=& \frac{1}{6}\left(1+\frac{1}{2}+\frac{1}{3}\right)-\frac{1}{6}=\frac{5}{36}.
\end{eqnarray}
\end{example}

Moreover, we can derive the Weighted Cumulative Entropy (WCE) of $X$ in terms of the second moment of the minimum order statistic in the following way
\begin{eqnarray}
\nonumber
\mathcal{CE}^w(X)&=&-\frac{x^2}{2}F(x)\log F(x)\big|_0^{+\infty}+\frac{1}{2}\int_0^{+\infty}x^2\log F(x)f(x)\mathrm dx+\\
\nonumber
& & \frac{1}{2}\int_0^{+\infty}x^2f(x)\mathrm dx \\
\nonumber
&=&\frac{1}{2}\int_0^{+\infty}x^2\log[1-(1-F(x))]f(x)\mathrm dx + \frac{1}{2}\mathbb E(X^2) \\
\nonumber
&=&-\frac{1}{2}\int_0^{+\infty}x^2\sum_{n=1}^{+\infty}\frac{(1-F(x))^n}{n}f(x)\mathrm dx+\frac{1}{2}\mathbb E(X^2) \\
\label{eq6}
&=&-\frac{1}{2}\sum_{n=1}^{+\infty}\frac{1}{n(n+1)}\mu^{(2)}_{1:n+1}+\frac{1}{2}\mathbb E(X^2),
\end{eqnarray}
where $\mu^{(2)}_{1:n+1}$ is the second moment of the smallest order statistic in a sample of size $n+1$, provided that $\lim_{x\to+\infty}-\frac{x^2}{2}F(x)\log F(x)=0$.

\section{Bounds}

Let us consider a sample with parent distribution $X$ such that $\mathbb E(X)=0$ and $\mathbb E(X^2)=1$. Hartley and David, 1954, and Gumbel, 1954, have shown that
\begin{equation}
\nonumber
\mu_{n:n}\leq\frac{n-1}{\sqrt{2n-1}}.
\end{equation}

We relate $\mu_{n:n}$ with the mean of the largest statistic order from the standard distribution. In fact, by normalizing the random variable $X$ with mean $\mu$ and variance $\sigma^2$ we get
\begin{equation}
\nonumber
Z=\frac{X-\mu}{\sigma}.
\end{equation}
Hence, the cdf $F_Z$ is given in terms of the cdf $F_X$ by
\begin{equation}
\nonumber
F_Z(x)=F_X(\sigma x+\mu).
\end{equation}
Then, cdf and pdf of the largest order statistic in a sample of size $n$ are
\begin{equation}
\nonumber
F_{Z_{n:n}}(x)=F_X^n(\sigma x+\mu), \ \ \ \ \ f_{Z_{n:n}}(x)=n F_X^{n-1}(\sigma x+\mu)f_X(\sigma x+\mu)\sigma.
\end{equation}
The mean of $X_{n:n}$ is given by
\begin{equation}
\nonumber
\mu_{n:n}=\mathbb E(X_{n:n})=n\int_0^{+\infty}xF_X^{n-1}(x)f_X(x)\mathrm dx.
\end{equation}
The mean of the largest statistic order from $Z$ is given by
\begin{eqnarray}
\nonumber
\mathbb E(Z_{n:n})&=&n\sigma\int_{-\frac{\mu}{\sigma}}^{+\infty}xF_X^{n-1}(\sigma x+\mu)f_X(\sigma x+\mu)\mathrm dx \\
\nonumber
&=&n\int_{0}^{+\infty}\frac{x-\mu}{\sigma}F_X^{n-1}(x)f_X(x)\mathrm dx \\
\nonumber
&=&\frac{n}{\sigma}\int_{0}^{+\infty}x F_X^{n-1}(x)f_X(x)\mathrm dx -\frac{n\mu}{\sigma}\int_{0}^{+\infty} F_X^{n-1}(x)f_X(x)\mathrm dx \\
\nonumber
&=&\frac{\mu_{n:n}-\mu}{\sigma}.
\end{eqnarray}

Using the Hartley-David-Gumbel bound for a non-negative parent distribution with mean $\mu$ and variance $\sigma^2$, we get
\begin{equation}
\label{eq7}
\mu_{n:n}=\sigma \mathbb E(Z_{n:n})+\mu \leq \sigma\frac{n-1}{\sqrt{2n-1}} +\mu.
\end{equation}

\begin{theorem}
\label{thm1}
Let $X$ be a non-negative random variable with mean $\mu$ and variance $\sigma^2$. Then, we obtain an upper bound for the CRE of $X$
\begin{equation}
\label{eq8}
\mathcal E(X)\leq \sum_{n=1}^{+\infty}\frac{\sigma}{(n+1)\sqrt{2n+1}}\simeq 1.21 \ \sigma.
\end{equation}
\end{theorem}

\proof
From \eqref{eq1} and \eqref{eq7} we get
\begin{eqnarray*}
\mathcal E(X)&=&\sum_{n=1}^{+\infty}\frac{1}{n(n+1)}\mu_{n+1:n+1}-\mathbb E(X) \\
&\leq&\sum_{n=1}^{+\infty}\frac{1}{n(n+1)}\left(\sigma\frac{n}{\sqrt{2n+1}} +\mathbb E(X)\right)-\mathbb E(X) \\
&=&\sum_{n=1}^{+\infty}\frac{\sigma}{(n+1)\sqrt{2n+1}}\simeq 1.21 \ \sigma,
\end{eqnarray*}
i.e., the upper bound given in \eqref{eq8}

\begin{remark}
\label{rem}
Since $X$ is non-negative we have that $\mu_{n+1:n+1}\ge0$, for all $n\in\mathbb N$. For this reason, using finite series approximations we get lower bounds for $\mathcal E(X)$:
\begin{equation}
\nonumber
\mathcal E(X)\ge \sum_{n=1}^{m}\frac{1}{n(n+1)}\mu_{n+1:n+1}-\mathbb E(X),
\end{equation}
for all $m\in\mathbb N$.
\end{remark}

\begin{remark}
\label{rem2}
Since $X$ is non-negative we have that $\mu_{1:n+1}\ge0$, for all $n\in\mathbb N$. For this reason, using finite series approximations we get upper bounds for $\mathcal {CE}(X)$:
\begin{equation}
\nonumber
\mathcal {CE}(X)\leq -\sum_{n=1}^{m}\frac{1}{n(n+1)}\mu_{1:n+1}+\mathbb E(X),
\end{equation}
for all $m\in\mathbb N$.
\end{remark}

\begin{theorem}
Let $X$ be DFR (decreasing failure rate). Then, we have the following lower bound for $\mathcal {CE}(X)$
\begin{equation}
\mathcal {CE}(X) \ge  \mathbb E(X)-\frac{\sqrt{\mathbb E(X^2)}}{\sqrt 2}\left(2-\frac{\pi^2}{6}\right).
\end{equation}
\end{theorem}

\proof
Let $X$ be DFR. From Theorem 12 of Rychlik (2001) we know that for a sample of size $n$, if 
\begin{equation}
\nonumber
\delta_j= \sum_{k=1}^j \frac{1}{n+1-k}\leq 2 \ \ \ j\in\{1,\dots,n\}
\end{equation}
then
\begin{equation}
\nonumber
\mathbb E(X_{j:n})\leq \frac{\delta_j}{\sqrt 2} \sqrt{\mathbb E (X^2)}.
\end{equation}
For $j=1$ we have $\delta_1=\frac{1}{n}\leq2$ for all $n\in\mathbb N$ and we get
\begin{equation}
\nonumber
\mathbb E(X_{1:n})\leq \frac{\sqrt{\mathbb E(X^2)}}{\sqrt 2 \ n}.
\end{equation}
Then, from \eqref{eq3} we get the following lower bound for $\mathcal {CE}(X)$
\begin{eqnarray}
\nonumber
\mathcal {CE}(X) &\ge& -\sum_{n=1}^{+\infty} \frac{1}{n(n+1)^2} \frac{\sqrt{\mathbb E(X^2)}}{\sqrt 2}+\mathbb E(X) \\
\nonumber
&=& \mathbb E(X)-\frac{\sqrt{\mathbb E(X^2)}}{\sqrt 2}\left(2-\frac{\pi^2}{6}\right).
\end{eqnarray}

\begin{remark}
We note that we can not provide an analogous bound for $\mathcal E(X)$ because $\delta_n\leq 2$ is not fulfilled for $n\ge 4$.
\end{remark}

David and Nagaraya, 2003, showed that if we have a sample $X_1,\dots,X_n$ with parent distribution $X$ symmetric about 0 with variance 1, then
\begin{equation}
\label{eq9}
\mathbb E(X_{n:n})\leq \frac{1}{2} \ nc(n),
\end{equation}
where
\begin{equation}
\nonumber
c(n)=\left[\frac{2\left(1-\frac{1}{\binom{2n-2}{n-1}}\right)}{2n-1}\right]^{\frac{1}{2}}.
\end{equation}

Using the bound \eqref{eq9} for a non-negative parent distribution symmetric about the mean $\mu$, with bounded support and variance $\sigma^2$, we get
\begin{equation}
\label{eq12}
\mu_{n:n}=\sigma \mathbb E(Z_{n:n})+\mu \leq \frac{1}{2}\sigma n c(n) +\mu.
\end{equation}

\begin{theorem}
Let $X$ be a symmetric non-negative random variable with bounded support, mean $\mu$ and variance $\sigma^2$. Then, we obtain an upper bound for the CRE of $X$
\begin{equation}
\label{eq13}
\mathcal E(X)\leq \frac{\sigma}{2}\sum_{n=1}^{+\infty}\frac{c(n)}{n}.
\end{equation}
\end{theorem}

\proof
From \eqref{eq1} and \eqref{eq12} we get
\begin{eqnarray*}
\mathcal E(X)&=&\sum_{n=1}^{+\infty}\frac{1}{n(n+1)}\mu_{n+1:n+1}-\mathbb E(X) \\
&\leq&\sum_{n=1}^{+\infty}\frac{1}{n(n+1)}\left(\frac{1}{2}\sigma (n+1) c(n+1) +\mathbb E(X)\right)-\mathbb E(X) \\
&=&\frac{\sigma}{2}\sum_{n=1}^{+\infty}\frac{c(n)}{n},
\end{eqnarray*}
i.e., the upper bound given in \eqref{eq13}.

About a symmetric distribution, Arnold and Balakrishnan, 1989, showed that if we have a sample $X_1,\dots,X_n$ with parent distribution $X$ symmetric about 0 with variance 1, then
\begin{equation}
\label{eq14}
\mathbb E(X_{n:n})\leq \frac{n}{\sqrt{2}} \sqrt{\frac{1}{2n-1}-B(n,n)},
\end{equation}
where $B(n,n)$ is the complete beta function.

Using the bound \eqref{eq14} for a non-negative parent distribution symmetric about the mean $\mu$ and with variance $\sigma^2$, we get
\begin{equation}
\label{eq15}
\mu_{n:n}=\sigma \mathbb E(Z_{n:n})+\mu \leq \sigma \frac{n}{\sqrt{2}} \sqrt{\frac{1}{2n-1}-B(n,n)} +\mu.
\end{equation}

\begin{theorem}
\label{tt}
Let $X$ be a symmetric non-negative random variable with mean $\mu$ and variance $\sigma^2$. Then, we obtain an upper bound for the CRE of $X$
\begin{equation}
\label{eq16}
\mathcal E(X)\leq \frac{\sigma}{\sqrt{2}}\sum_{n=1}^{+\infty}\frac{1}{n} \sqrt{\frac{1}{2n+1}-B(n+1,n+1)}.
\end{equation}
\end{theorem}

\proof
From \eqref{eq1} and \eqref{eq15} we get
\begin{eqnarray*}
\mathcal E(X)&=&\sum_{n=1}^{+\infty}\frac{1}{n(n+1)}\mu_{n+1:n+1}-\mathbb E(X) \\
&\leq&\sum_{n=1}^{+\infty}\frac{1}{n(n+1)}\left(\sigma \frac{n+1}{\sqrt{2}} \sqrt{\frac{1}{2n+1}-B(n+1,n+1)} +\mathbb E(X)\right)-\mathbb E(X) \\
&=&\frac{\sigma}{\sqrt{2}}\sum_{n=1}^{+\infty}\frac{1}{n} \sqrt{\frac{1}{2n+1}-B(n+1,n+1)},
\end{eqnarray*}
i.e., the upper bound given in \eqref{eq16}.

\begin{example}
Let us consider a sample with parent distribution $X\sim N(0,1)$. From Harter, 1961, we get the values of the mean of the largest order statistic for samples of size less than 100. Hence, we compare the finite series approximation of \eqref{eq1} and \eqref{eq16} and we expect the same relation given in Theorem \ref{tt} because truncated terms are negligible. We get the following result
\begin{eqnarray}
\nonumber
\mathcal E(X)&\simeq& \sum_{n=1}^{99}\frac{1}{n(n+1)}\mu_{n+1:n+1}\simeq 0.87486 \\
\nonumber
&<& \frac{1}{\sqrt{2}}\sum_{n=1}^{99}\frac{1}{n} \sqrt{\frac{1}{2n+1}-B(n+1,n+1)}\simeq 0.94050.
\end{eqnarray}
\end{example}

From \eqref{eq1} and \eqref{eq3} we get the following expression for the sum of the cumulative residual entropy and the cumulative entropy
\begin{equation}
\label{eq10}
\mathcal{E}(X)+\mathcal{CE}(X)= \sum_{n=1}^{+\infty}\frac{1}{n(n+1)}(\mu_{n+1:n+1}-\mu_{1:n+1}).
\end{equation}
Calì et al., 2017, showed a connection among \eqref{eq10} and the partition entropy studied by Bowden, 2007.

\begin{theorem}
\label{thm4}
We have the following bound for the sum of the CRE and the CE
\begin{equation}
\mathcal{E}(X)+\mathcal{CE}(X)\leq  \sum_{n=1}^{+\infty}\frac{\sqrt{2}\ \sigma }{n\sqrt{n+1}}\simeq 3.09\ \sigma.
\end{equation}
\end{theorem}

\proof
From Theorem 3.24 of Arnold and Balakrishnan, 1989, we know the following bound for the difference between the expectation of the largest and the smallest order statistics from a sample of size $n+1$
\begin{equation}
\label{eq11}
\mu_{n+1:n+1}-\mu_{1:n+1} \leq \sigma \sqrt{2(n+1)},
\end{equation}
and so using \eqref{eq11} in \eqref{eq10} we get the following bound for the sum of the CRE and the CE
\begin{equation}
\nonumber
\mathcal{E}(X)+\mathcal{CE}(X)\leq  \sum_{n=1}^{+\infty}\frac{\sigma \sqrt{2(n+1)}}{n(n+1)}=\sum_{n=1}^{+\infty}\frac{\sqrt{2}\ \sigma }{n\sqrt{n+1}}\simeq 3.09\ \sigma.
\end{equation}

About a symmetric distribution, Arnold and Balakrishnan, 1989, showed that if we have a sample $X_1,\dots,X_n$ with parent distribution $X$ symmetric about the mean $\mu$ with variance 1, then
\begin{equation}
\label{eq17}
\mathbb E(X_{n:n})-\mathbb E(X_{1:n})\leq n \sqrt{2} \sqrt{\frac{1}{2n-1}-B(n,n)},
\end{equation}
where $B(n,n)$ is the complete beta function.

Using the bound \eqref{eq17} for a non-negative parent distribution symmetric about the mean $\mu$ and with variance $\sigma^2$, we get
\begin{equation}
\label{eq18}
\mu_{n:n}-\mu_{1:n}=\sigma \left(\mathbb E(Z_{n:n})-\mathbb E(Z_{1:n})\right) \leq \sigma n\sqrt{2} \sqrt{\frac{1}{2n-1}-B(n,n)}.
\end{equation}

\begin{theorem}
Let $X$ be a symmetric non-negative random variable with mean $\mu$ and variance $\sigma^2$. Then, we obtain an upper bound for the sum of the CRE and the CE of $X$
\begin{equation}
\label{eq19}
\mathcal{E}(X)+\mathcal{CE}(X)\leq \sqrt{2} \ \sigma\sum_{n=1}^{+\infty}\frac{1}{n} \sqrt{\frac{1}{2n+1}-B(n+1,n+1)}.
\end{equation}
\end{theorem}

\proof
From \eqref{eq10} and \eqref{eq18} we get
\begin{eqnarray*}
\mathcal{E}(X)+\mathcal{CE}(X)&=& \sum_{n=1}^{+\infty}\frac{1}{n(n+1)}(\mu_{n+1:n+1}-\mu_{1:n+1}) \\
&\leq&\sum_{n=1}^{+\infty}\frac{1}{n(n+1)}\left(\sigma (n+1)\sqrt{2} \sqrt{\frac{1}{2n+1}-B(n+1,n+1)}\right) \\
&=&\sqrt{2} \ \sigma \sum_{n=1}^{+\infty}\frac{1}{n}\sqrt{\frac{1}{2n+1}-B(n+1,n+1)},
\end{eqnarray*}
i.e., the upper bound given in \eqref{eq19}.

In Table \ref{table1} we present some applications of the bounds obtained in this section to important distributions in the reliability theory.

\begin{table}
\caption{Some bounds for known distributions.}
\centering
%% \tablesize{} %% You can specify the fontsize here, e.g., \tablesize{\footnotesize}. If commented out \small will be used.
\begin{tabular}{cccccc}
\toprule
Support  & CDF &	$\mathcal{E}(X)$  & Bound thm.\ref{thm1} & $\mathcal{E}(X)+\mathcal{CE}(X)$ & Bound thm.\ref{thm4}  \\
\midrule
$x>0$ & $1-\exp(-\lambda x)$ &$\frac{1}{\lambda}$ & $\frac{1.21}{\lambda}$ & $\frac{\pi^2}{6 \ \lambda}$ & $\frac{3.09}{\lambda}$ \\
$0<x<a$ & $\frac{x}{a}$ & $\frac{a}{4}$ &$ \frac{1.21 \ a}{2\sqrt{3}}$ & $\frac{a}{2}$ & $\frac{3.09 \ a}{2\sqrt{3}}$ \\
$x\in (0,1)$ & $\frac{1}{x^2}\exp\left(2\left(1-\frac{1}{x}\right)\right)$ & $0.1549$ &$ 0.1999$ & $0.2936$ & $0.5105$ \\
$x>0$ & $1-\frac{1}{(x+1)^3}$ & $0.75$ & $1.0479$ & $1.1115$ & $2.6759$ \\
$x\in (0,1)$ & $x^2$ & $0.1869$ &$ 0.2852$ & $0.4091$ & $0.7283$ \\
$x\in (0,+\infty)$ & $\exp\left(-\frac{1}{\exp(x)-1}\right)$ & $0.9283$ &$ 1.1238$ & $1.5246$ & $2.87$ \\
\bottomrule
\label{table1}
\end{tabular}
\end{table}

\begin{acknowledgements}
Francesco Buono and Maria Longobardi are partially supported by the GNAMPA research group of INdAM (Istituto Nazionale di Alta Matematica) and MIUR-PRIN 2017, Project "Stochastic Models for Complex Systems" (No. 2017 JFFHSH).
\end{acknowledgements}

\section*{Conflict of interest}

The authors declare that they have no conflict of interest.


\begin{thebibliography}{99}
\footnotesize

\bibitem{arnold}
Arnold, B. C., Balakrishnan, N. (1989). {\em Relations, Bounds and Approximations for Order Statistics.} Springer, New York, NY. 

\bibitem{bowden}
Bowden, R. (2007). Information, measure shifts and distribution diagnostics. \emph{Statistics}, \textbf{46}(2), 249--262. 

\bibitem{cali}
Calì, C., Longobardi, M., Ahmadi, J. (2017). Some properties of cumulative Tsallis entropy. \emph{Physica A}, \textbf{486}, 1012--1021. 

\bibitem{CaLoNa}
Calì, C.; Longobardi, M.; Navarro, J. (2020). Properties for generalized cumulative past measures of information. 
\emph{Probab. Eng. Inform. Sci.} \textbf{34}, 92--111.

\bibitem{CaLoPsa}
Calì, C.; Longobardi, M.; Psarrakos, G. (2019). A family of weighted distributions based on the mean inactivity time and cumulative past entropies.
\emph{Ricerche Mat.}, (in press),  
 doi:10.1007/s11587-019-00475-7.

\bibitem{david}
David, H. A., Nagaraya, H. N. (2003). {\em Order Statistics.} John Wiley \& Sons Inc. 

\bibitem{dicrescenzo}
Di Crescenzo, A., Longobardi, M. (2009). On cumulative entropies. \emph{Journal of Statistical Planning and Inference}, \textbf{139}, 4072--4087.

\bibitem{gumbel}
Gumbel, E. J. (1954). The maxima of the mean largest value and of the range. \emph{The Annals of Mathematical Statistics}, \textbf{25}, 76--84.

\bibitem{harter}
Harter, H. L. (1961). Expected Values of Normal Order Statistics. \emph{Biometrika}, \textbf{48}(1), 151--165. 

\bibitem{hartley}
Hartley, H. O., David, H. A.  (1954). Universal bounds for mean range and extreme observation. \emph{The Annals of Mathematical Statistics}, \textbf{25}, 85--99.

\bibitem{Lo2014}
Longobardi, M. (2014). Cumulative measures of information and stochastic orders.
\emph{Ricerche Mat.}  \textbf{63}, 209--223.

\bibitem{mirali}
Mirali, M., Baratpour, S., Fakoor, V., (2016). On weighted cumulative residual entropy. \emph{Communications in Statistics - Theory and Methods}, \textbf{46}(6), 2857--2869.

\bibitem{mirali2}
Mirali, M., Baratpour, S., (2017). Some results on weighted cumulative entropy. \emph{Journal of the Iranian Statistical Society}, \textbf{16}(2), 21--32.

\bibitem{rao}
Rao, M., Chen, Y., Vemuri, B., Wang, F. (2004). Cumulative Residual Entropy: A New Measure of Information. \emph{IEEE Transactions on Information Theory}, \textbf{50}, 1220--1228.

\bibitem{rychlik}
Rychlik, T. (2001). {\em Projecting Statistical Functionals.} Springer, New York, NY. 

\bibitem{shannon}
Shannon, C. E. (1948). A mathematical theory of communication. \emph{Bell System Technical Journal}, \textbf{27}, 379--423. 


\end{thebibliography}
\end{document}